\documentclass[11pt,a4paper]{article}

\usepackage{amsmath}
\usepackage{amsthm}
\usepackage{amssymb}

\usepackage[latin1]{inputenc}

\usepackage[a4paper=true,%
breaklinks=true,%
colorlinks=true,%
linkcolor={blue},%
citecolor={red},%
pdftitle={Spacetime causality in the study of the Hankel transform},%
pdfauthor={Jean-Francois Burnol},%
pdfkeywords={Klein-Gordon equation;%
 Hankel and Fourier transforms; Scattering},%
pdfsubject={Mathematical physics},%
pdfstartview=FitH,
]{hyperref}

\def\From{From}

\newtheorem{theorem}{Theorem}

\theoremstyle{definition}

\newcommand{\RR}{{\mathbb R}}
\newcommand{\CC}{{\mathbb C}}
\newcommand{\cH}{{\mathcal H}}
\newcommand{\cA}{{\mathcal A}}
\newcommand{\cB}{{\mathcal B}}
\newcommand{\cS}{{\mathcal S}}
\newcommand{\cK}{{\mathcal K}}
\DeclareMathOperator{\Un}{{\bf 1}}

\let\wh=\widehat

\let\cal\mathcal

\newcommand{\ddx}{\frac{d}{dx}}

\newcommand{\ddu}{\frac{d}{du}}
\newcommand{\dpt}{\frac{\partial}{\partial t}}
\newcommand{\dpx}{\frac{\partial}{\partial x}}

\newcommand{\dphidt}{\frac{\partial\phi}{\partial t}}
\newcommand{\dphidx}{\frac{\partial\phi}{\partial x}}
\newcommand{\dphidu}{\frac{\partial\phi}{\partial u}}
\newcommand{\dphidv}{\frac{\partial\phi}{\partial v}}
\newcommand{\dpsidt}{\frac{\partial\psi}{\partial t}}
\newcommand{\dpsidx}{\frac{\partial\psi}{\partial x}}
\newcommand{\dpsidu}{\frac{\partial\psi}{\partial u}}

\begin{document}

\title{Spacetime causality in the study of the Hankel transform}

\author{Jean-Fran\c cois Burnol}

\date{\ }

\maketitle
\medskip

\begin{abstract}
We study Hilbert space aspects of the Klein-Gordon equation
in two-dimensional spacetime. We associate to its
restriction to a spacelike wedge a scattering from the past
light cone to the future light cone, which is then shown to
be (essentially) the Hankel transform of order zero. We
apply this to give a novel proof, solely based on the
causality of this spatio-temporal wave propagation, of the
theorem of de~Branges and V.~Rovnyak concerning Hankel pairs
with a support property. We recover their isometric
expansion as an application of Riemann's general method for
solving Cauchy-Goursat problems of hyperbolic type.
\end{abstract}

\medskip

\begin{footnotesize}
\begin{quote}
keywords:
  Klein-Gordon equation; Hankel and Fourier transforms; Scattering.

\medskip

Universit\'e Lille 1\\ 
UFR de Math\'ematiques\\ 
Cit\'e scientifique M2\\ 
F-59655 Villeneuve d'Ascq\\ 
France\\
burnol@math.univ-lille1.fr\\
\bigskip
v1: 25 Sept. 2005; Final: 5 March 2006\\
\end{quote}
\end{footnotesize}

\medskip

\setlength{\normalbaselineskip}{18pt}
\baselineskip\normalbaselineskip
\setlength{\parskip}{12pt}

\section{Introduction}

We work in two-dimensional spacetime with metric $c^2dt^2 -
dx^2$. We shall use units such
  that $c=1$. Points are denoted $P = (t,x)$. And the
  d'Alembertian operator $\square$ is
  $\frac{\partial^2}{\partial t^2} -
  \frac{\partial^2}{\partial x^2}$. We consider the
  Klein-Gordon equation (with $m=1$, $\hbar = 1$; actually
  we shall only study the classical wave field, no
  quantization is involved in this paper):
\begin{equation}
\square \phi + \phi = 0 \end{equation}
We have an energy density:\footnote{as this paper is
  principally of a mathematical nature, we do not worry
  about an overall $\frac12$ factor.}
\begin{equation} \cal E = |\phi|^2 + \left|\dphidx\right|^2 +
  \left|\dphidt\right|^2\end{equation}
which gives a conserved quantity:
\begin{equation} E = \frac1{2\pi} \int_{-\infty}^{+\infty} \cal
  E(\phi)(t,x)\;dx\;,\end{equation} 
in the sense that if the Cauchy data at time $t=0$ has
  $E<\infty$ then $E$ is finite (and constant\dots) at all
  times (past and future). We shall mainly work with such
  finite energy solutions.
Although we failed in locating a reference for the following
  basic observation, we can not imagine it to be novel:

  \begin{theorem}\label{theo:1}
    If $\phi$ is a finite energy solution to the
  Klein-Gordon equation then:
 \[ \lim_{t\to+\infty} \int_{|x|>t} \cal
  E(\phi)(t,x)\;dx = 0\;.\]
  \end{theorem}

Obviously this would be completely wrong for the zero mass
  equation. We shall give a (simple) self-contained proof,
  because it is the starting point of all that we do
  here. Let us nevertheless state that the result follows
  immediately from Hörmander's fine pointwise estimates
  (\cite{hor1, hor2}; see also the paper of S.~Klainerman
  \cite{klain} and the older papers of S.~Nelson \cite{nel1,
  nel2}.)  I shall not reproduce the strong pointwise
  results of Hörmander, as they require notations and
  preliminaries.  Let me simply mention that Hörmander's
  Theorem 2.1 from \cite{hor1} can be applied to the
  positive and negative frequency parts of a solution with
  Cauchy data which is gaussian times polynomial. So theorem
  \ref{theo:1} holds for them, and it holds then in general,
  by an approximation argument.

The energy conservation follows from:
\begin{equation} \dpt \cal E  + \dpx \cal P = 0
\qquad\text{with}\qquad \cal P = -  \dphidx \overline{\dphidt} -
  \overline{\dphidx} \dphidt \end{equation}
If we apply Gauss' theorem to the triangle with vertices $O =
  (0,0)$, $A = (t,t)$, $B = (t,-t)$, we obtain ($t>0$):
\[ \int_{|x|\leq t} \cal E(\phi)(t,x) dx = \int_{-t}^0 (|\phi(|x|,x)|^2 +
  \left| \ddx \phi(|x|,x)\right|^2) \,dx + \int_{0}^t (|\phi(x,x)|^2 +
  \left| \ddx \phi(x,x)\right|^2) \,dx\]
This proves that   $\int_{|x|>t} \cal
  E(\phi)(t,x)\;dx$ decreases as $t\to+\infty$. It shows
  also that
  theorem \ref{theo:1} is equivalent to:
\begin{equation}\label{eq:123}
 E = \frac1{2\pi} \int_{-\infty}^0 (|\phi(|x|,x)|^2 +
  \left| \ddx \phi(|x|,x)\right|^2 ) \,dx + \frac1{2\pi}
  \int_{0}^\infty ( |\phi(x,x)|^2 + 
  \left| \ddx \phi(x,x)\right|^2 ) \,dx
\end{equation}

Otherwise stated, there is a unitary representation of
  $\phi$ on the future light cone.
Here
  is now the basic idea: as
  solutions to hyperbolic equations propagate causally,
  equation \eqref{eq:123}
  gives a unitary representation from the Hilbert space of
  Cauchy data at time $t=0$ with support in $x\geq0$ to the
  Hilbert space of functions $p(v) = \phi(v,v)$ on
  $[0,+\infty[$ with squared norm $ \frac1{2\pi} \int_{0}^\infty (
  |p(v)|^2 + |p'(v)|^2 ) \,dv$. Instead of Cauchy data
  vanishing for $x<0$, it will be useful to use Cauchy data
  invariant under $(t,x)\to(-t,-x)$. Then $p$ will be
  considered as an even, and $p'$ as an odd, function, and $
  \frac1{2\pi} \int_{0}^\infty ( |p(v)|^2 + |p'(v)|^2 )
  \,dv$ will be $\frac12 E(\phi)$, for $\phi(t,x) =
  \phi(-t,-x)$.
We can also consider the
  past values $g(u) = \phi(-u,u)$, $t = -u$, $x =u$, $0\leq
  u <\infty$. So there is a unitary map from such $g$'s to
  the $p$'s:

  \begin{theorem}\label{theo:2}
    Let $g(u)$, $u>0$, and $p(v)$, $v>0$ be such that
    $\int_0^\infty
  |g(u)|^2+|g'(u)|^2 du < \infty$, $\int_0^\infty
  |p(v)|^2+|p'(v)|^2 dv < \infty$. The necessary and
  sufficient condition for $A(r) = \sqrt{r}\,
  g(\frac{r^2}2)$ and $B(s) = - \sqrt{s}\, p'(\frac{s^2}2)$
  to be Hankel transforms of order zero of one another
  ($A(r) = \int_0^\infty \sqrt{rs} J_0(rs) B(s) \,ds$) is
  for $g$ and $p$ to be the values on the past and future
  boundaries of the Rindler wedge $0<|t|<x$ of a finite
  energy solution $\phi(t,x)$ of the Klein-Gordon equation
  ($g(u) = \phi(-u,u)$, $p(v) = \phi(v,v)$.) For any $a>0$
  the vanishing on $0<x<2a$ of the Cauchy data for
  $\phi(t,x)$ at $t=0$ is the necessary and sufficient
  condition for the simultaneous vanishing of $g(u)$ for
  $0<u<a$ and of $p(v)$ for $0<v< a$.
  \end{theorem}

The statements relative to the support
  properties are corollaries to the relativistic
  causality of the propagation of solutions to the
  Klein-Gordon equation. Regarding the function $B$, 
if $k(v)=-p'(v)$ vanishes identically on $(0,a)$, then $p(v)$ is constant
there, and this constant has to
  be $0$ if $g(u)$ is also identically zero on $(0, a)$: indeed
  the finite energy solution $\phi$ is continuous on
  spacetime (this follows from the well-known explicit
  formulas \eqref{eq:classic}).
We employed temporarily $A(r) = \sqrt{r}\, g(\frac{r^2}2)$ and $B(s) = -
  \sqrt{s}\, p'(\frac{s^2}2)$ in the statement of Theorem \ref{theo:2} in
  order to express the matter with the zero order Hankel
  transform. It proves more natural  to stay with $g(u)$ and
  $k(v)=-p'(v)$. They are connected by the integral formula:
$g(u) = \int_0^\infty J_0(2\sqrt{uv}) k(v)\,dv$, so this
  motivates the definition of the $\cH$ transform:
  \begin{equation}
    \label{eq:3}
    \cH(f)(x) = \int_0^\infty J_0(2\sqrt{xy}) f(y)\,dy
  \end{equation}
The $\cH$ transform is a unitary operator on
  $L^2(0,+\infty;dx)$ which is self-reciprocal. As is
  well-known $\sqrt x e^{-\frac12 x^2}$ is an invariant
  function for the Hankel transform of order zero, so, for
  the $\cH$ transform we have $e^{-x}$ as invariant function
  in $L^2(0,\infty;dx)$. The $\cH$ operator is
  ``scale-reversing'': by this we mean that $\cH(f(\lambda
  y))(x) = \lambda^{-1} \cH(f)(\lambda^{-1} x)$, or,
  equivalently, that the operator $\cH\cdot I$ is scale
  invariant, where $I$ is the unitary operator
  $f(x)\mapsto\frac1x f(\frac 1x)$. As we explain later,
  $\cH$ is the unique scale-reversing operator on
  $L^2(0,\infty;dx)$ having among its self-reciprocal
  functions the function $e^{-x}$. Let us restate Theorem
  \ref{theo:2} as it applies to $\cH$:

  \begin{theorem}\label{theo:2b}
    Let $\phi(t,x)$ be a finite energy solution of the
  Klein-Gordon equation. Let $g(u) = \phi(-u,u)$ for $u>0$
  and $p(v)= \phi(v,v)$ for $v>0$ be the values taken by
  $\phi$ on the past, respectively future, boundaries of the
  Rindler wedge $0<|t|<x$. Then $k(v) = -p'(v)$ is the $\cH$
  transform of $g(u)$: $k(v) = \int_0^\infty J_0(2\sqrt{uv})
  g(u)\,du$. For any $a>0$ the vanishing for $0<x<2a$,
  $t=0$, of the Cauchy data for $\phi(t,x)$ is the necessary
  and sufficient condition for the simultaneous vanishing of
  $g(u)$ for $0<u<a$ and $p(v)$ for $0<v<a$.
  \end{theorem}

In this manner a link has been established between the
  relativistic causality and a mathematical theorem of
  de~Branges \cite{bra64}, and V.~Rovnyak \cite{rov1} (see
  further \cite{rov2}). They proved an explicit isometric
  representation of $L^2(0,+\infty; dx)$ onto
  $L^2(0,+\infty; dy)\oplus L^2(0,+\infty; dy)$,
  $h\mapsto(f,g)$, such that the zero order Hankel transform
  on $L^2(0,+\infty;dx)$ is conjugated to the simple map
  $(f,g)\to (g,f)$, and such that the pair $(f(y),g(y))$
  vanishes identically on $(0,a)$ if and only $h(x)$ and its
  Hankel transform of order zero both identically vanish on
  $(0,a)$. Their formulas ((5) and (7) of \cite{bra64}
  should be corrected to read as
      (3) and (2) of \cite{rov1}) are:
      \begin{subequations}
\begin{align}
f(y) &= \int_y^\infty h(x) J_0(y\sqrt{x^2 -
  y^2})\sqrt{xy}\,dx\\
g(y) &= h(y) - \int_y^\infty h(x)y\,\frac{J_1(y\sqrt{x^2
  -y^2})}{\sqrt{x^2 - y^2}}\,\sqrt{xy}\,dx\\
h(x) &= g(x) + \int_0^x f(y) J_0(y\sqrt{x^2 -
  y^2})\sqrt{xy}\,dy
 - \int_0^x g(y)\frac{J_1(y\sqrt{x^2 -
  y^2})}{\sqrt{x^2 - y^2}} \sqrt{xy}\,y\,dy\\
\int_0^\infty &|h(x)|^2\,dx = \int_0^\infty (|f(y)|^2 +
  |g(y)|^2)\,dy
\end{align}
      \end{subequations}

We shall give an independent, self-contained proof, that
  these formulas are mutually compatible and have the stated
  relation to the Hankel tranform of order zero.  The main
  underlying idea has been to realize
  the Hankel transform of order zero as a scattering related
  to a causal propagation of waves. The support condition
  initially considered by de~Branges and Rovnyak has turned out to be
  related to relativistic
causality, and the looked-after scattering has been realized as
  the transition  from the past to the future boundary of the Rindler
  wedge $0<|t|<x$.  Also, in the technique of proof we
  apply, in a perhaps
unusual manner, the classical Riemann method
(\cite[IV\textsection 1]{john},
  \cite[VI\textsection5]{courant}) from the theory of
  hyperbolic equations. Let us reformulate here the
  isometric
expansion of de~Branges-Rovnyak into a version which applies
  to the $\cH$ transform. For this we write, for $x>0$,
\[ h(x) = \sqrt x\; k(\frac{x^2}2),\qquad f(x) = \sqrt x
  F(x^2), \qquad g(x) = \sqrt x G(x^2) \]
Then the equations above become:
\begin{subequations}
\begin{align}\label{eq:100a}
    F(x) &= 
  \int_{x/2}^\infty J_0(\sqrt{x(2v-x)})k(v)\,dv\\
\label{eq:100b}
G(x) &= k(\frac{x}2) -
  \int_{x/2}^\infty x\,
  \frac{J_1(\sqrt{x(2v-x)})}{\sqrt{x(2v-x)}}\,
  k(v)\,dv\\
\label{eq:100c}
k(v)  &= G(2v) + \frac12\int_0^{2v} 
  J_0(\sqrt{x(2v - x)}) F(x) \,dx  - \frac12\int_0^{2v}
  x\, \frac{J_1(\sqrt{x(2v - x)})}{\sqrt{x(2v - x)}} G(x)
  \,dx\\
\label{eq:100d}
\int_0^\infty 2&|k(v)|^2\,dv = \int_0^\infty (|F(x)|^2 +
  |G(x)|^2)\,dx 
\end{align}
\end{subequations}
The de~Branges Rovnyak theorem is thus the equivalence between
  equations \eqref{eq:100a}, \eqref{eq:100b} and
  \eqref{eq:100c}, the validity of \eqref{eq:100d}, the
  fact that the pair $(F,G)$ is identically zero on $(0,2a)$
  if and only if both $k$ and $\cH(k)$ vanish identically on
  $(0,a)$, and finally the fact that permuting $F$ and $G$ is equivalent to
  $k\leftrightarrow \cH(k)$.

It proves convenient to work with the first order
``Dirac'' system:
\begin{subequations}
\begin{align}
\dpsidt - \dpsidx = + \phi\\
\dphidt + \dphidx = - \psi
\end{align}
\end{subequations}
Let us write $\left[\begin{smallmatrix} \psi(0,x)\\
  \phi(0,x)\end{smallmatrix}\right] = \left[\begin{smallmatrix} G(x)\\
  F(x)\end{smallmatrix}\right]$. We shall use $K(\psi,\phi) = \frac1{2\pi}
\int_{-\infty}^{+\infty}
  (|F(x)|^2 + |G(x)|^2)\,dx$ as the Hilbert space (squared)
  norm. We shall require $\dphidx$ and $\dpsidx$ to be
  in $L^2$ at $t=0$ (then $\phi$ and $\psi$ are continuous
  on space-time). Our previous $E(\phi)$ is not invariant
  under Lorentz boosts:
  it is only the first component of a Lorentz vector
  $(E(\phi),P(\phi))$ (see equation \eqref{eq:P} for the
  expression of $P$). And it turns out that in fact $K(\psi,\phi) =
  E(\phi) - P(\phi) = E(\psi) + P(\psi)$. The point
is that in order to define an action of the
  Lorentz group on the solutions of the Dirac system it is
  necessary to rescale in opposite ways $\psi$
  and $\phi$. When done symmetrically, $K$ then becomes an
  invariant under the Lorentz boosts. This relativistic
  covariance of the spinorial quantity
  $\left[\begin{smallmatrix} \psi \\
  \phi \end{smallmatrix}\right]$ is important for the proof
  of the next theorem:

  \begin{subequations}
\begin{theorem}\label{theo:4}
Let $F$ and $G$ be two functions with $\int_0^\infty
  |F|^2+|F'|^2+|G|^2+|G'|^2 dx < \infty$. 
  Let $\left[{\psi\atop \phi}\right]$ be the unique solution
  in the Rindler wedge $x>|t|>0$ of the
  first order system:
\begin{align}
\dpsidt - \dpsidx &= + \phi\\
\dphidt + \dphidx &= - \psi
\end{align}
with Cauchy data $\phi(0,x) = F(x)$, $\psi(0,x) =
  G(x)$. The boundary values:
\[ g(u) = \phi(-u,u) \quad (u>0),\qquad\text{and}\qquad k(v) =
  \psi(v,v) \quad (v>0),\]
verify $\int_0^\infty |g(u)|^2 + |g'(u)|^2 du
  <\infty$, $\int_0^\infty |k(v)|^2 + |k'(v)|^2 dv < \infty$ and are a
  $\cH$ transform pair.  For any $a>0$ the identical
  vanishing of $F(x)$ and $G(x)$ for $0<x<2a$ is equivalent to the
  identical vanishing of $g(u)$ for $0< u< a$ and of $k(v)$ for
  $0< v <a$. All $\cH$ pairs with $\int_0^\infty |g|^2 +
  |g'|^2 du <\infty$, $\int_0^\infty |k|^2 + |k'|^2 dv < \infty$ are
  obtained in this way.  The functions $F(x)$, $G(x)$,
$g(u)$ and
  $k(v)$ are related by the following formulas:
\begin{align}
    F(x) &=
  \int_{x/2}^\infty J_0(\sqrt{x(2v-x)})k(v)\,dv
= g(\frac{x}2) -
  \int_{x/2}^\infty x\,
  \frac{J_1(\sqrt{x(2u-x)})}{\sqrt{x(2u-x)}}\, g(u)\,du\\
G(x) &= k(\frac{x}2) -
  \int_{x/2}^\infty x\, 
  \frac{J_1(\sqrt{x(2v-x)})}{\sqrt{x(2v-x)}}\,k(v)\,dv
=
  \int_{x/2}^\infty J_0(\sqrt{x(2u-x)})g(u)\,du\\
g(u) &= F(2u) + \frac12\int_0^{2u}
  J_0(\sqrt{x(2u - x)}) G(x) \,dx - \frac12\int_0^{2u} x\,
  \frac{J_1(\sqrt{x(2u - x)})}{\sqrt{x(2u - x)}} F(x) \,dx\\
k(v) &= G(2v) + \frac12\int_0^{2v}
  J_0(\sqrt{x(2v - x)}) F(x) \,dx - \frac12\int_0^{2v} x\,
  \frac{J_1(\sqrt{x(2v - x)})}{\sqrt{x(2v - x)}} G(x) \,dx\\
\int_0^\infty 2&|k(v)|^2\,dv = \int_0^\infty (|F(x)|^2 +
  |G(x)|^2)\,dx = \int_0^{\infty} 2|g(u)|^2\,du\\
k(v) &= \int_0^\infty J_0(2\sqrt{uv})g(u)\,du\qquad\qquad g(u) =
\int_0^\infty J_0(2\sqrt{uv})k(v)\,dv
\end{align}
The integrals converge as improper Riemann integrals. 
\end{theorem}
  \end{subequations}

The Lorentz boost parameter can
serve as ``time'' as $K$ is conserved under it. In this
  manner  going-over from $\phi$ on the past light cone
  to $\psi$ on the future light cone becomes a
  scattering. We shall explain its formulation in the
  Lax-Phillips \cite{laxphillips} terminology.  

In conclusion we can say that
  this paper identifies the unique scale reversing operator
  $\cH$ on $L^2(0,+\infty;\,dx)$ such that $e^{-x}$ is
  self-reciprocal as the scattering from the past (positive
  $x$)-light-cone to the future (positive $x$)-light-cone
  for finite energy solutions of the Dirac-Klein-Gordon
  equation in two-dimensional space-time. Some further
  observations and remarks will be found in the concluding
  section of the paper. The operator $\cH$, which is
  involved in some functional equations of number theory, is
  studied further by the author in \cite{hankel}.

\section{Plane waves}

Throughout this paper we shall use the following light cone
  coordinates, which are positive on the right wedge:
  \begin{subequations}
\begin{gather}
 v= \frac{x+t}2\quad\quad u = \frac{x-t}2 \\
x = u + v\qquad t = - u + v\qquad t^2 - x^2 = 4 (-u) v \qquad
\square = - \frac{\partial^2}{\partial u \partial v}
\end{gather}
  \end{subequations}
We write sometimes $\phi(t,x) = \phi[u,v]$. 

Let us begin the proof of Theorem \ref{theo:1}. We can build
a solution to the Klein-Gordon equation by
  superposition of plane waves:
  \begin{subequations}
\begin{gather}
 \phi(t,x) = \int_{-\infty}^{+\infty}  e^{+i(\lambda u - \frac1\lambda
  v)}\alpha(\lambda)\,d\lambda = \int_{-\infty}^{+\infty}  e^{- i (\omega t - \mu
  x)}\alpha(\lambda)\,d\lambda\\
 \text{with}\qquad \omega =
  \frac12(\lambda+\frac1\lambda),\qquad \mu =
  \frac12(\lambda - \frac1\lambda)
\end{gather} 
  \end{subequations}
The full range
  $-\infty<\lambda<+\infty$ allows to keep track
  simultaneously of the ``positive frequency'' ($\lambda>0$,
  $\omega\geq1$), and ``negative frequency'' ($\lambda<0$,
  $\omega\leq-1$) parts.

At first we only take $\alpha$ to be a smooth, compactly
  supported function of $\lambda$, vanishing identically in
  a neighborhood of $\lambda=0$. Then the corresponding
  $\phi$ is a smooth, finite energy solution of the Klein
  Gordon equation. Let us compute this energy. At $t=0$ we
  have 
\[ \phi(0,x) = \int_{-\infty}^{+\infty}  e^{+i \mu x} \alpha(\lambda)\,d\lambda
  \qquad \dphidt(0,x) = - i \int_{-\infty}^{+\infty}  e^{+i\mu x}
  \frac12(\lambda+\frac1\lambda) \alpha(\lambda)\,d\lambda\]
So we will apply Plancherel's theorem, after the change of
  variable $\lambda\to \mu$. We must be careful that if
  $\lambda$ is sent to $\mu$, then $\lambda' =
  -\frac1\lambda$, is too. Let $\lambda_1 >0$ and $\lambda_2
  <0$ be the ones being sent to $\mu$. Let us also define:
\[ a(\mu) =
  \frac{\alpha(\lambda_1)}{\frac12(1+\frac1{\lambda_1^2})},
  \qquad
 b(\mu) =
  \frac{\alpha(\lambda_2)}{\frac12(1+\frac1{\lambda_2^2})}\]
Then:
\[ \phi(0,x) = \int_{-\infty}^{+\infty}  e^{+i \mu x} ( a(\mu) + b(\mu) ) d\mu
\qquad  \dphidt(0,x) = - i \int_{-\infty}^{+\infty}  e^{+i\mu x}
  \frac12(\lambda_1+\frac1{\lambda_1}) ( a(\mu) - b(\mu) )d\mu\]
\[ \frac1{2\pi} \int_{-\infty}^{+\infty}  (|\phi|^2  + |\dpx \phi|^2) dx = \int_{-\infty}^{+\infty}  | a(\mu) + b(\mu)
  |^2 (1 + \mu^2)\, d\mu \]
\[ \frac1{2\pi} \int_{-\infty}^{+\infty}  |\dpt \phi|^2 dx = \int_{-\infty}^{+\infty}  | a(\mu) - b(\mu)
  |^2 \left(\frac12(\lambda_1+\frac1{\lambda_1})\right)^2\,
  d\mu \]
Observing that $1+\mu^2 =
  \left(\frac12(\lambda_1+\frac1{\lambda_1})\right)^2 =
  \left(\frac12(\lambda_2+\frac1{\lambda_2})\right)^2  $, this gives
\[ E(\phi) = 2 \int_{-\infty}^{+\infty}  (|a(\mu)|^2 + |b(\mu)|^2) \lambda_1^2
  \left(\frac12(1+\frac1{\lambda_1^2})\right)^2 d\mu\]
\[ = 2\int_0^\infty |a(\mu)|^2 \lambda_1^2
  \left(\frac12(1+\frac1{\lambda_1^2})\right)^3 \,d\lambda_1
 + 2\int_{-\infty}^0 |b(\mu)|^2 \lambda_2^2
  \left(\frac12(1+\frac1{\lambda_2^2})\right)^3\,d\lambda_2\]
\[ = 2\int_0^\infty |\alpha(\lambda)|^2 \frac12 (1+\lambda^2) \, d\lambda
  +2\int_{-\infty}^0 |\alpha(\lambda)|^2 \frac12 (1+\lambda^2)\, d\lambda\] 
\begin{equation} E(\phi) = \int_{-\infty}^{+\infty} |\alpha(\lambda)|^2 (1 +
  \lambda^2)\,d\lambda\end{equation}

Let us now compute the energy on the future light cone. We write
  $g(u) = \phi(-u,u)$, $p(v) = \phi(+v,v)$. We
  have:
  \begin{equation}
    \label{eq:5}
    g(u) = \int_{-\infty}^{+\infty}  e^{+ i \lambda u}
    \alpha(\lambda)\,d\lambda 
  \end{equation}
Let $\alpha = \alpha_+ + \alpha_-$ be the decomposition of
  $\alpha$ as the sum of $\alpha_+$,
  belonging to the
  Hardy space of the upper half-plane $\Im(\lambda)>0$ and
 of $\alpha_-$, belonging to the Hardy space of the
  lower half-plane. We have:
  \begin{subequations}
    \begin{equation}\label{eq:gandgprime}
      \frac1{2\pi} \int_{-\infty}^0 |g(u)|^2 \,du = \int_{-\infty}^{+\infty} 
  |\alpha_+(\lambda)|^2 \,d\lambda
 \end{equation}
We must be careful regarding $\frac1{2\pi} \int_{-\infty}^0
 |g'(u)|^2 \,du$. We have 
$g'(u)\Un_{u<0}(u) - g(0)\delta(u) =
 \int_{-\infty}^{+\infty}  e^{+ i \lambda u} 
    i\lambda \alpha_+(\lambda)\,d\lambda$ as a distribution identity
so
\[ \frac1{2\pi} \int_{-\infty}^0
 |g'(u)|^2 \,du = \int_{-\infty}^{+\infty} \left|
   \lambda\alpha_+(\lambda) +
   \frac{g(0)}{2\pi i}\right|^2\,d\lambda\]
On the other
 hand $(t\alpha(t))_+(\lambda) = \frac1{2\pi
   i}\int_{-\infty}^\infty \frac{t\alpha(t)\,dt}{t-\lambda}
   =  \lambda\alpha_+(\lambda) + \frac{g(0)}{2\pi i}$ so the
   formula is:
   \begin{equation}
     \label{eq:11}
     \frac1{2\pi} \int_{-\infty}^0
 |g'(u)|^2 \,du = \int_{-\infty}^{+\infty} 
  |(t\alpha)_+(\lambda)|^2 \,d\lambda 
   \end{equation}
  \end{subequations}
Similarly, as $p(v) = \int_{-\infty}^{+\infty}  e^{+ i \lambda v}
  \alpha(-\frac1\lambda)\frac1{\lambda^2}\,d\lambda$, and
  defining  $\beta(\lambda) = \alpha(-
  \frac1\lambda)\frac1{\lambda^2}$ we obtain:
\[ \frac1{2\pi} \int_0^\infty  (|p(v)|^2 + |p'(v)|^2) dv =
  \int_{-\infty}^{+\infty}  
  (|\beta_-(\lambda)|^2  + |(t\beta)_-(\lambda)|^2)\,d\lambda\]
Now, from $t\beta(t) = \frac1t \alpha(\frac{-1}t)$ it is
  seen that $(t\beta)_-(\lambda) = \frac1\lambda
  \alpha_-(\frac{-1}\lambda)$, so $\int_{-\infty}^{+\infty}
  |(t\beta)_-(\lambda)|^2\,d\lambda = \int_{-\infty}^{+\infty}
  |\alpha_-(\lambda)|^2\,d\lambda$. And, as $t\alpha(t) =
  \frac1t \beta(\frac{-1}t)$ one has in a similar manner
  $\int_{-\infty}^{+\infty} 
  |(t\alpha)_+(\lambda)|^2\,d\lambda = \int_{-\infty}^{+\infty}
  |\beta_+(\lambda)|^2\,d\lambda$.
Combining, we get 
  \begin{align*}
    &\frac1{2\pi} \int_{-\infty}^0 (|g(u)|^2  +
  |g'(u)|^2) \,du + \frac1{2\pi} \int_0^\infty  (|p(v)|^2 +
  |p'(v)|^2) dv\\
 &= \int_{-\infty}^{+\infty} 
   (|\alpha_+(\lambda)|^2 + |\beta_+(\lambda)|^2 +
  |\beta_-(\lambda)|^2 +|\alpha_-(\lambda)|^2)  \,d\lambda
\end{align*}
and
  as the two Hardy spaces are mutually perpendicular in
  $L^2(-\infty,+\infty;d\lambda)$ we finally obtain:
\[ \int_{-\infty}^{+\infty} |\alpha(\lambda)|^2 (1+ \lambda^2)
  \,d\lambda\]
as the energy on the future light cone.

So, with this, the theorem that $E(\phi)$ is entirely on the
  future light cone is proven for the $\phi$'s corresponding
  to $\alpha$'s which are smooth and compactly supported
  away from $\lambda = 0$. Obviously the Cauchy data for
  such $\phi$'s is a dense subspace of the full initial data
  Hilbert space. As energy is conserved as $t\to\infty$, the
  fact that $\lim_{t\to\infty} \int_{|x|>t} \cal E(\phi) dx
  =0$ holds for all finite energy $\phi$'s then follows by
  approximation.  Furthermore we see that a finite energy
  solution is uniquely written as a wave packet:
\begin{equation} \phi(t,x) = \int_{-\infty}^{+\infty}
  e^{+i(\lambda u - \frac1\lambda v)} 
  \alpha(\lambda)\,d\lambda\qquad E(\phi) =  \int_{-\infty}^{+\infty}  (1+
  \lambda^2)|\alpha(\lambda)|^2\,d\lambda < \infty
\end{equation}
At this stage Theorem \ref{theo:1} is established.

When studying the Klein-Gordon equation in the right wedge
$x>0$, $|t|<x$, we can arbitrarily extend the Cauchy data to
$x<0$. Setting it to $0$ on $x<0$, however, will be
compatible with the finite energy condition only if
$\phi(0,0^+) = 0$. If this is the case then this choice
makes $g(u)$ vanish for $u<0$ and $p(v)$ vanish for
$v<0$, which is a Hardy space constraint on $\alpha$, in
fact it means that $\alpha$ belongs to the Hardy space of
the lower half-plane.  Another manner to extend the Cauchy
data to $x<0$ is to make it invariant under the $PT$
operation $(t,x)\to (-t,-x)$. This has the advantage, if
$\int_0^\infty |\phi(0,x)|^2 + |\frac\partial{\partial
x}\phi(0,x)|^2\,dx<\infty$, to produce Cauchy data of finite
energy on the full line $-\infty<x<\infty$. The condition on
$\alpha$ is to be even. We shall often use this
convention when studying the Klein-Gordon equation in the
right wedge.

\section{Energy and momentum}

The momentum density $\cal P = -  \dphidx \overline{\dphidt} -
  \overline{\dphidx} \dphidt$ also satisfies a conservation
  law:
\[ \dpt \cal P + \dpx \left( - |\phi|^2 + \left|\dphidx\right|^2 +
  \left|\dphidt\right|^2 \right) = 0 \]
So 
\begin{equation}\label{eq:P}
 P = - \frac1{2\pi} \int_{-\infty}^{+\infty} \left( \dphidx
\overline{\dphidt} +
  \overline{\dphidx} \dphidt\right) dx 
\end{equation}
is also a conserved
  quantity. We have:
  \begin{subequations}
\begin{align}\label{eq:emoinsp} E - P &= \frac1{2\pi}
  \int_{-\infty}^{+\infty} \left ( |\phi|^2 + 
  \left|\dphidx + \dphidt \right|^2 \right) \,dx\\
  E + P &= \frac1{2\pi} \int_{-\infty}^{+\infty} \left ( |\phi|^2 +
  \left|\dphidx - \dphidt \right|^2 \right) \,dx
\end{align}
  \end{subequations}
Applying Gauss' theorem to $\cal P$ we obtain for $t>0$:
\[ \int_{|x|\leq t} \cal P(\phi)(t,x) dx = \int_{-t}^0 (- |\phi(-x,x)|^2 +
  \left| \ddx \phi(-x,x)\right|^2) \,dx + \int_{0}^t (|\phi(x,x)|^2 -
  \left| \ddx \phi(x,x)\right|^2) \,dx\]
The integral of $| \cal P |$ for $|x|>t$ tends to zero for
  $t\to+\infty$ as it is bounded above by the one for $\cal
  E$. So:
  \begin{equation}
    \label{eq:7}
    P = \frac1{2\pi} \int_{-\infty}^0 ( - |g(u)|^2 +
  |g'(u)|^2)\,du + \frac1{2\pi} \int_{0}^\infty  ( |p(v)|^2 -
  |p'(v)|^2)\,dv
  \end{equation}
with, again, $g(u) = \phi(-u,u)$, $p(v) =
  \phi(+v,v)$. Hence:
  \begin{subequations}
\begin{align}\label{eq:1}
 E - P &= \frac1{\pi} \int_{-\infty}^0 
  |g(u)|^2\,du +  \frac1{\pi} \int_{0}^\infty |p'(v)|^2\,dv\\
 E + P &= \frac1{\pi} \int_{-\infty}^0 
  |g'(u)|^2\,du +  \frac1{\pi} \int_{0}^\infty |p(v)|^2\,dv 
\end{align}
  \end{subequations}
\From\ \eqref{eq:gandgprime} and the similar formulas relative
  to $p$ we can express all four integrals in terms of
  $\alpha(\lambda)$. Doing so we find after elementary steps:
  \begin{equation}
    \label{eq:9}
E - P = 2 \int_{-\infty}^{+\infty}  |\alpha(\lambda)|^2\,d\lambda\qquad E + P
  = 2  \int_{-\infty}^{+\infty}  \lambda^2\,|\alpha(\lambda)|^2\,d\lambda
  \end{equation}
So:
\begin{equation} P = \int_{-\infty}^{+\infty} (\lambda^2 -
  1)|\alpha(\lambda)|^2\,d\lambda\end{equation}
This confirms that a $\lambda$ with $|\lambda|\geq1$
gives a ``right-moving''
  component of the wave packet (its  phase is constant for
  $\omega t - \mu x = C$, $\omega =
  \frac12(\lambda+\frac1\lambda)$, $\mu = \frac12(\lambda-
  \frac1\lambda)$.) The values of $\lambda$ with
$|\lambda|\leq1$ give ``left-moving'' wave components. As a check,
  we can observe that it is impossible to have a purely
  right-moving packet with vanishing Cauchy data for $t=0$,
  $x<0$, because as we saw above, for such Cauchy data
  $\alpha$ has to belong to the Hardy space of the lower
  half-plane and can thus (by a theorem of Wiener) not
  vanish identically on $(-1,1)$. A purely right-moving
  packet starting entirely on $x>0$ would have a hard time
  hitting the light cone, and this would emperil  Theorem
  \ref{theo:1}. Such wave-packets exist for the zero-mass
  equation, one way of reading Theorem \ref{theo:1} is to
  say that they don't exist for non-vanishing real mass.

Let us consider the effect of a Lorentz boost on $E$ and
  $P$. We take $\Lambda = e^\xi$ ($\xi\in\RR$) and replace
  $\phi$ by:
\begin{subequations}
\begin{gather}
 \phi_\Lambda(t,x) = \phi( \cosh(\xi) t +
  \sinh(\xi) x, \sinh(\xi) t + \cosh(\xi) x )\\
\phi_\Lambda[u,v] = \phi[\frac1\Lambda u, \Lambda v]\\
 g_\Lambda(u) = \phi_\Lambda[u,0] = g(\frac1\Lambda
  u)\qquad p_\Lambda(v) = p(\Lambda v)\label{eq:gscale}\\
 \alpha(\lambda) \mapsto \alpha_\Lambda(\lambda) = \Lambda
  \alpha(\Lambda \lambda) \\
 E_\Lambda - P_\Lambda = \Lambda\cdot(E -
  P)\qquad\qquad  E_\Lambda + P_\Lambda = \frac1\Lambda (E +
  P)\\
 E_\Lambda = \cosh(\xi) E - \sinh(\xi) P\\
 P_\Lambda = -\sinh(\xi) E + \cosh(\xi) P
\end{gather}
\end{subequations}

So the conserved quantities $E$ and $P$ are not Lorentz
  invariant but the Einstein rest mass squared $E^2 - P^2$
  is.

\section{Scale reversing operators}

We begin the proof of Theorem \ref{theo:2}.
Let us consider the manner in which the function $g(u)$
  for $u>0$ is related to the function $p(v)>0$. We know
  that they are in unitary correspondence for the norms
  $\int_{u>0} |g|^2+|g'|^2 \,du$ and $\int_{v>0}
  |p|^2+|p'|^2\,dv$, and the formulas \eqref{eq:1} for $E-P$
  and $E+P$ suggest that one should pair $g$ with $p'$ and
  $g'$ with $p$. In fact if we take into consideration the
   wave which has values $\phi(t,x) = e^{-|x|}$
  for space-like points, we are rather led to pair $g$ with
  $-p'$ and $g'$ with $-p$ (the values of $\phi$ at
  time-like points are more involved and we don't need to
  know about them here; suffice it to say that certainly
  $e^{-x}$ solves Klein-Gordon, so it gives the unique
  solution in the right wedge with $\phi(0,x) = e^{-x}$,
  $\dphidt(0,x) = 0$.)

Let us denote by $\cH$ the operator which acts as $g\mapsto
-p'$,
  on even $g$'s. Under a Lorentz boost:
 $g\mapsto g_\Lambda(u) = g(\frac1\Lambda
  u)$, $ -p'\mapsto -\Lambda p'(\Lambda v) $
and also the assignment $g\mapsto -p'$ is unitary for the
$L^2$ norm:
\[ g(u) = \int_{-\infty}^{+\infty}  e^{iu\lambda}
  \alpha(\lambda)\,d\lambda\qquad  
 p(v) = \int_{-\infty}^{+\infty}  e^{i \lambda v}
  \alpha(-\frac1\lambda)\frac1{\lambda^2}\,d\lambda 
  \]
\[ -p'(v) = -i \int_{-\infty}^{+\infty}  e^{i \lambda v}
  \alpha(-\frac1\lambda)\frac1{\lambda}\,d\lambda \]
Going from $g$ to $\alpha$ is unitary, from $\alpha$ to
  $-i\alpha(-\frac1\lambda)\frac1{\lambda}$ also, and back
  to $-p'$ also, in the various $L^2$ norms. So the assignment from
  $g$ to $-p'$ is unitary.

Identifiying the $L^2$ space on $u>0$ with
  the $L^2$ space on $v>0$, through $v=u$, $\cH$ is a
  unitary operator on $L^2(0,+\infty;\,du)$. Furthermore it
  is ``scale reversing'': we say that an operator $\cK$
  (bounded, more generally, closed) is scale reversing if
  its composition $\cK I$ with $I: g(u)\mapsto
  \frac{g(1/u)}u$ commutes with the unitary group of scale
  changes $g\mapsto \sqrt \Lambda g(\Lambda u)$. The Mellin
  transform $g\mapsto\wh g(s) = \int_0^\infty g(u)
  u^{-s}\,du$, for $s= \frac12+i\tau$, $\tau\in\RR$, is the
  additive Fourier transform of $e^{t/2}g(e^{t})\in
  L^2(-\infty,+\infty;dt)$.  The operator $\cK I$
  commutes with multiplicative translations hence is
  diagonalized by the Mellin transform: we have a certain
  (bounded for $\cK$ bounded) measurable function $\chi$
  on the critical line $\Re(s)= \frac12$ such that for any
  $g(u) \in L^2(0,\infty;du)$, and almost everywhere on the
  critical line:
\[ (\cK g)^\wedge(s) = (\cK I(Ig))^\wedge(s) = \chi(s)
  (Ig)^\wedge(s) = \chi(s) \wh g(1-s) \]
Let us imagine for a minute that we know a $g$ which is
  invariant under $\cK$ and which, furthermore has $\wh
  g(s)$ almost everywhere non vanishing (by a theorem of
  Wiener, this means exactly that the linear span of its orbit under the
  unitary group of scale changes is dense in $L^2$). Then we
  know $\chi(s)$ hence, we know $\cK$. So $\cK$ is
  uniquely determined by the knowledge of one such invariant
  function. 

In the case of our operator $\cH$ which goes from the
  data of $g(u)$, $u>0$, to the data of $k(v) = -p'(v)$,
  $v>0$, where $g$ and $p$ are the boundary values of a
  finite energy solution of the Klein-Gordon equation in the
  right wedge, we know that it is indeed unitary, scale
  reversing, and has $e^{-u}$ as a self-reciprocal function
  (so, here, $\chi(s) = \frac{\Gamma(1-s)}{\Gamma(s)}$).

On the other hand the Hankel transform of order zero is
  unitary, scale reversing, and has $\sqrt u e^{-u^2/2}$ as
  self-reciprocal invariant function. So we find that the
  assignment of $-\sqrt v\;k'(\frac{v^2}2)$ to $\sqrt u\;
  g(\frac{u^2}2)$ is exactly the Hankel transform of order
  zero. This may also be proven directly by the method we
  will employ in section \ref{sec:riemann}.

\section{Causality and support conditions}

The Theorem \ref{theo:2} is almost entirely proven: if the
  Cauchy data vanishes identically for $0<x<2a$, then by
  unicity and causal propagation, $g(u) = \phi(-u,u)$
  vanishes identically for $0<u<a$ and $p(v) = \phi(+v,v)$
  vanishes identically for $0<v<a$. Conversely, if $A$ and
  $B$ from Theorem \ref{theo:2} vanish identically for
  $0<r,s<\sqrt{2a}$, then $g(u)$ and $-p'(v)$ vanish
  identically for $0<u< a$ and $0<v<a$. We explained in the
  introduction that $p$ itself also vanishes identically for
  $0<v< a$. Then $\phi[u,v] = \iint_{0\leq r\leq u\atop
  0\leq s\leq v} \phi[r,s]\,drds$ for $0\leq u\leq
 a$, $0\leq v\leq a$, hence $\phi$ vanishes
  identically in this range, and the Cauchy data for $\phi$
  at $t=0$, $0<x<2a$, vanishes identically. The proof of
  Theorem \ref{theo:2} (hence also in its equivalent form
  \ref{theo:2b}) is complete.

We would like also to relax the finite energy condition on
  $\phi$. Let us imagine that our $g$, say even, is only
  supposed $L^2$. It has an $L^2$ Fourier transform
  $\alpha$ such that $g(u) = \int_{-\infty}^{+\infty}  e^{+i u\lambda}
  \alpha(\lambda)\,d\lambda$. Let us approximate $\alpha$ by an $L^2$
  converging sequence of
  $\alpha_n$'s, corresponding to finite
  energy Klein-Gordon solutions $\phi_n$. We have by
  \eqref{eq:emoinsp} and \eqref{eq:9}:
\[ \frac1{2\pi} \int_{-\infty}^{+\infty} \left (
  |\phi_n - \phi_m|^2 +
  \left|\frac{\partial(\phi_n-\phi_m)}{\partial x} + 
\frac{\partial(\phi_n-\phi_m)}{\partial t}\right|^2
  \right)\,dx
= 2\int_{-\infty}^{+\infty}  |\alpha_n - \alpha_m|^2\,d\lambda \]
So the $\phi_n$ converge for $t=0$ in the $L^2$ sense, and
  also the $\frac{\partial\phi_n}{\partial x} +
\frac{\partial\phi_n}{\partial t}$. We can then consider, as
  is known to exist, the distribution solution $\phi$ with
  this Cauchy data.

Let us suppose that we start from an even $g$ which,
  together with its $\cH$ transform, vanish in
  $(0,a)$. First we show that we can find, with
  $0<b_n<1$, $b_n\to1$, a sequence of $g_n$'s, such that $g_n'$
  is in $L^2$, and $g_n\to g$ in $L^2$, with the $g_n$'s
  satisfying  the
  support condition for $(0,b_n a)$. We obtain such
  $g_n$ by multiplicative convolution of $g$ with a test
  function supported in $(b_n,\frac 1{b_n})$. At the level of
  Mellin transforms, this multiplies by a Schwartz
  function. As $u\ddu$ corresponds to multiplication by $-s$
  certainly the $u\ddu$ of our $g_n$'s are in $L^2$. But
  then $\ddu g_n $ itself is in $L^2$ as we know that it vanishes
  in $(0,b_n a)$. And its $\cH$ transform also
  vanishes there.

So the corresponding $\phi_n$'s will have for $t=0$ 
  vanishing Cauchy data in intervals only arbitrarily slightly smaller
  than $(0,2a)$. The $L^2$ functions
  $\phi(0,x)$ and $(\dphidx + \dphidt)(0,x)$ will thus vanish
  identically, in $(0,2a)$. Conversely if we have two $L^2$
  functions $L$ and $M$ vanishing in $(0,2a)$ we can
  approximate then by Schwartz functions $L_n$ and $M_n$
  vanishing in $(0,b_n 2a)$ ($0<b_n<1$, $b_n\to1$), solve
  the Cauchy problem with data $\phi = L_n$ and $\dphidx +
  \dphidt = M_n$ at $t=0$, consider the corresponding
  $g_n$'s which vanish identically for $0<u<b_n a$ and get
  an $L^2$ limit $g$ vanishing identically in $(0,a)$. The
  $\cH$ transform of $g$ will be the limit in $L^2$ of the
  $\cH$ transforms of the $g_n$, so it will also vanish in
  $(0,a)$.

In conclusion the space-time representation of Hankel pairs
  with  support condition as given in Theorem \ref{theo:2}
  extends to the general case of $L^2$ Hankel pairs if one
  allows Klein-Gordon solutions of possibly infinite energy
  but such that  $\phi(0,x)$ and $\dphidx(0,x) +
  \dphidt(0,x)$ are in $L^2$.

\section{The Dirac system and its associated scattering}

We return to finite energy
  solutions which are associated to functions $\alpha$
  verifying the condition $\int_{-\infty}^{+\infty}
  (1+\lambda^2)|\alpha(\lambda)|^2\,d\lambda <\infty$. Let
  us consider in fact a pair $\left[\begin{smallmatrix}
  \psi\\ \phi\end{smallmatrix}\right]$ of such finite energy
  solutions
satisfying the first order system:

\begin{subequations}\label{eq:dirac}
\begin{align}
\dpsidt - \dpsidx = + \phi\qquad\qquad \dpsidu = -\phi\\
\dphidt + \dphidx = - \psi\qquad\qquad \dphidv = -\psi
\end{align}
\end{subequations}

If $\alpha$ corresponds to $\phi$ and $\beta$ corresponds to
  $\psi$, then there is the relation:
$ \alpha(\lambda) = -i\lambda \beta(\lambda)$
so we must have $\int_{-\infty}^{+\infty}  \frac1{\lambda^2}
  |\alpha(\lambda)|^2\,d\lambda <\infty$. To enact  a Lorentz
  boost we could imagine replacing 
  $\phi$ and $\psi$ by 
\[ \phi( \cosh(\xi) t +
  \sinh(\xi) x, \sinh(\xi) t + \cosh(\xi) x ) =
\phi[e^{-\xi} u, e^\xi v]\]
\[ \psi( \cosh(\xi) t +
  \sinh(\xi) x, \sinh(\xi) t + \cosh(\xi) x ) = \psi[e^{-\xi} u, e^\xi v]
\]
but this does not give a  solution of the Dirac type
system \eqref{eq:dirac}. To obtain a solution we must
rescale $\phi$, or
  $\psi$, or both. We choose:\footnote{this
  conflicts with our previous notation
  $\phi_\Lambda[u,v]=\phi[\frac1\Lambda u,\Lambda v]$; no
  confusion should arise.} 
\begin{equation}\label{eq:spinor}
 \phi_\xi[u,v] =  e^{-\xi/2}\; \phi[e^{-\xi} u, e^\xi
v]\qquad
\psi_\xi[u,v] = e^{\xi/2} \; \psi[e^{-\xi} u, e^\xi
v]
\end{equation}
In other words, if we want to consider our $\phi$ as a
component of such a system we must cease treating it as a
scalar. It is a (spinorial) quantity which transforms as
indicated under a Lorentz boost. We note further that with
this modification both $E(\phi) - P(\phi)$ and 
 $E(\psi) + P(\psi)$ are Lorentz invariant. In fact
they are identical: $E(\phi) - P(\phi) = \frac1{2\pi}
\int_{-\infty}^{+\infty}  |\phi|^2 +
  \left|\dphidx + \dphidt \right|^2 
\,dx$, $E(\psi) + P(\psi) = \frac1{2\pi}
\int_{-\infty}^{+\infty}  |\psi|^2 +
  \left|\dpsidx - \dpsidt \right|^2 \,dx$, hence:
\begin{equation} E(\phi) - P(\phi)= E(\psi) + P(\psi) = \frac1{2\pi}
  \int_{-\infty}^{+\infty} \left(|\phi(0,x)|^2 + |\psi(0,x)|^2
  \right)\; dx\end{equation}

We again focus on what happens in the right wedge. Thus, we can as
  well take $\phi$ to be  $PT$ invariant. But then as $\psi
  = -\dphidv$, $\psi$ must acquire a sign under the $PT$
  transformation: $\psi(-x,-t) = - \psi(x,t)$. So the
  function $g(u)   = \phi(-u,u) = \phi[u,0]$ is even but the
  function $k(v) = \psi(v,v) = \psi[0,v]$ is odd. In fact
  $k(v) = -p'(v)$ with our former notation. So we know that
  the $PT$ invariant $\phi$ is uniquely determined by $g(u)$
  for $u>0$ which gives under the $\cH$ transform the
  function $k(v)$ for $v>0$ which must be considered odd and
  correspond to the $PT$ anti-invariant $\psi$. We note that
  if $k(0^+)\neq0$ then this $\psi$ is not of finite
  energy. Using only that $\phi$ is finite energy, we have 
from equation \eqref{eq:1}:
\[E(\phi)  - P(\phi) = \frac1{\pi} \int_{-\infty}^0 
  |g(u)|^2\,du +  \frac1{\pi} \int_{0}^\infty |k(v)|^2\,dv\]
\[ \frac1{2\pi} \int_{-\infty}^{+\infty} \left(|\phi(0,x)|^2 + |\psi(0,x)|^2
  \right)\; dx = E(\phi) - P(\phi) = \frac1{\pi} \int_0^{\infty} 
  |g(u)|^2\,du +  \frac1{\pi} \int_{0}^\infty |k(v)|^2\,dv\]
  \begin{subequations}
    \begin{align}\label{eq:13a}
   \int_0^\infty \left(|\phi(0,x)|^2 + |\psi(0,x)|^2
  \right)\; dx &= 2  \int_0^{\infty} 
  |g(u)|^2\,du \\ \label{eq:13b}
   \int_0^\infty \left(|\phi(0,x)|^2 + |\psi(0,x)|^2
  \right)\; dx &=  2  \int_{0}^\infty |k(v)|^2\,dv      
    \end{align}
  \end{subequations}

We now begin the proof of Theorem \ref{theo:4}. To prove
that $\int_0^\infty |F(x)|^2+|G(x)|^2\,dx = 2\int_0^\infty
|g(u)|^2\,du = 2\int_0^\infty |k(v)|^2\,dv$, we extend $F$
to be even and $G$ to be odd. Then $\phi$ is $PT$ even of
finite energy, and $\psi$ is $PT$ odd and equations
\eqref{eq:13a} and \eqref{eq:13b} apply. Note that if
$G(0^+)\neq0$ then $\psi$ is not of finite energy but only
the fact that $\phi$ is of finite energy was used for
\eqref{eq:13a} and \eqref{eq:13b}.  That $k=\cH(g)$ and
$\int_0^\infty |g(u)|^2
+ |g'(u)|^2 du <\infty$ hold are among our previous
results. If we choose $G$ to be even and $F$ to be odd, then
it is $\psi$ which is of finite energy and so $\int_0^\infty
|k(v)|^2 + |k'(v)|^2 dv < \infty$ holds true. We can also
prove $\int_0^\infty |g|^2
+ |g'|^2 du <\infty$, $\int_0^\infty |k|^2 + |k'|^2 dv <
\infty$ after extending $F$ and $G$ such that
$\int_{-\infty}^\infty
  |F|^2+|F'|^2+|G|^2+|G'|^2 dx < \infty$ so that both
$\phi$ and $\psi$ are then of finite energy. The boundary
values $g(u)$, $u>0$, and $k(v)$, $v>0$ do not depend on
choices.  Furthermore the vanishing of $F$ and $G$ on
$(0,2a)$ at $t=0$ is equivalent by our previous arguments to
the vanishing of $g$ and $k$ on $(0,a)$. To show that all
$\cH$ pairs with $\int_0^\infty |g|^2 + |g'|^2 du <\infty$,
$\int_0^\infty |k|^2 + |k'|^2 dv < \infty$ are obtained, let
$k_1$ be the odd function with $k_1(v) = k(v) -
k(0^+)e^{-v}$ for $v>0$ and let $g_1$ be the even function
with $g_1(u) = g(u) - k(0^+)e^{-u}$ for $u\geq 0$. Then $k_1
= \cH(g_1)$ and $\int_{-\infty}^\infty |g_1|^2 + |g_1'|^2 du
<\infty$ and $\int_{-\infty}^\infty |k_1|^2 + |k_1'|^2 dv <
\infty$. They thus correspond to $\phi_1$ and $\psi_1$ both
of finite energy. We define for $x>0$: $F(x) = \phi_1(0,x) +
k(0^+) e^{-x}$ and $G(x) = \psi_1(0,x) + k(0^+) e^{-x}$, it
then holds that $\int_{0}^\infty
  |F|^2+|F'|^2+|G|^2+|G'|^2 dx < \infty$ and
$\left[\begin{smallmatrix} \psi \\ \phi
\end{smallmatrix}\right] = \left[\begin{smallmatrix} \psi_1
+ k(0^+) e^{-x}\\
 \phi_1 + k(0^+)e^{-x} \end{smallmatrix}\right]$ is
the unique solution in the Rindler wedge of the Dirac system
with Cauchy data $\left[\begin{smallmatrix} G \\ F
\end{smallmatrix}\right]$ on $x>0$, $t=0$, and it has $g(u)$
and $k(v)$ as boundary values. To complete the proof of
Theorem \ref{theo:4} there only remains to show the formulas
relating $F$, $G$, $g$, and $k$ and this will be done in the
next section.

On the Hilbert space $L^2(0,\infty;\,dx)\oplus
  L^2(0,\infty;\,dx)$ of the pairs $(F,G)$, we can define
  a unitary group $U(\xi)$, $-\infty<\xi<\infty$, as
  follows: we define its action at first for $(F,G)$ with
  $F', G' \in L^2$. Let $\left[{\psi\atop
  \phi}\right]$ be the solution of first order system
  \eqref{eq:dirac}  such that  $\phi(0,x) = F(x)$, $\psi(0,x) =
  G(x)$. Then we take:
\begin{equation}
U(\xi)(F,G) = (\left.\phi_{-\xi}\right|_{t=0}, \left.\psi_{-\xi}\right|_{t=0})
\end{equation}
where \eqref{eq:spinor} has been used. As
  $\xi$ increases from $-\infty$ to $+\infty$ this has the
  effect of transporting $\phi$ and $\psi$ forward along the
  Lorentz boosts trajectories. We can also implement
  $U(\xi)$ as a unitary group acting on the $L^2$ space of
  the $g(u) = \phi(-u,u)$ functions, or on the space of the
  $k(v) = \psi(v,v)$ functions. We then have,
taking into account \eqref{eq:spinor} (and $-\xi$):
  \begin{equation}
    g_\xi(u) = e^{\frac\xi2} g(e^{\xi} u)\qquad\qquad 
    k_\xi(v) = e^{-\frac\xi2} k(e^{-\xi} v)
  \end{equation}

Following the terminology of Lax-Phillips
  \cite{laxphillips} (the change of variable $u\to \log(u)$
  would reduce to the additive language of
  \cite{laxphillips}) we shall say that $(F,G)\mapsto I(g)$
  provides an incoming (multiplicative) translation
  representation ($U(\xi)$ moves the graph of
  $e^{y/2}I(g)(e^{y}) = e^{-y/2}g(e^{-y})$ to the right by
  an amount of additive time $\xi$) and $(F,G)\mapsto k$ is
  an outgoing translation representation. We use $(Ig)(u) =
  \frac1u g(\frac1u)$ as it is translated by $U(\xi)$ in the
  same direction as $k$. The assignment $Ig\to k$ will be
  called the ``scattering matrix'' $\cal S$ (it is canonical
  only up to a translation in ``time'', which means here
  only up to a scale change in $u$).  With our previous
  notation it is $S = \cH I$. Let us give a ``spectral''
  representation of $S$. For this we represent $g$ as a
  superposition of (multiplicative) harmonics, $g(u) =
  \frac1{2\pi} \int_{\Re(s)=\frac12} \wh g(s)
  u^{s-1}\,|ds|$, with $\wh g(s) = \int_0^\infty g(u)\,
  u^{-s}\,du$, $s=\frac12+i \tau$. Then the unitary operator
  $S$ will be represented as multiplication by a unit
  modulus function $\chi(s)$. Multiplication by $\chi(s)$
  must send the Mellin transform $\Gamma(s)$ of $I(e^{-u})$
  to the Mellin transform $\Gamma(1-s)$ of $e^{-u}$, in
  other words:
  \begin{equation}\label{eq:chi}
    \chi(s) = \frac{\Gamma(1-s)}{\Gamma(s)}
  \end{equation}
We thus see that the first order system in the wedge of two
  dimensional space-time provides an interpretation of this
  function (for $\Re(s) = \frac12$) as a scattering matrix.
  To obtain the Hankel transform of order zero, and not its
  succ\'edan\'e $\cH$, one writes $s = \frac14 + \frac
  w2$, where again $\Re(w) = \frac12$. In fact, with our
  normalizations, the scattering matrix corresponding to the
  tranform $g(t)\mapsto f(u)=\int_0^\infty
  \sqrt{ut}J_0(ut)g(t)\,dt$ is
 the function $2^{\frac12 -w}\frac{\Gamma(\frac34 - \frac
  w2)}{\Gamma(\frac14+\frac w2)}$ on the critical line
  $\Re(w) = \frac12$.

\section{Application of Riemann's method}\label{sec:riemann}

The completion of the proof of Theorem \ref{theo:4} will now
  be provided.
I need to briefly review
Riemann's method (\cite[IV\textsection 1]{john},
  \cite[VI\textsection5]{courant}), although it is such a
  classical thing, as I will use it in a special manner
  later. In the case of the (self-adjoint) Klein-Gordon
  equation $\frac{\partial^2 \phi}{\partial u\partial v} =
  +\phi$, $t^2 - x^2 = 4(-u)v$, Riemann's method combines:
  \begin{enumerate}
  \item whenever $\phi$ and $\psi$ are two solutions,  the
  differential form 
$ \omega = \phi\dpsidu\,du + \psi\dphidv\,dv$
is
  closed,
\item it is advantageous to use either for $\phi$ or for
  $\psi$ the special solution (Riemann's function) $R(P,Q)$
which reduces to the constant value $1$ on
  each of characteristics issued from a given point
  $P$. Here $R(P,Q) = R(P-Q,0) = R(Q-P,0)$, $R((t,x),0) =
  J_0(\sqrt{t^2 - x^2}) = J_0(2\sqrt{-uv})$.
  \end{enumerate}

Usually one uses Riemann's method to solve for $\phi$ when
  its Cauchy data is given on a 
  curve transversal to the characteristics. But one can
  also use it when the data is on the characteristics
  (Goursat problem). Also,
  one usually symmetrizes the formulas obtained in combining
  the information from using $\phi\frac{\partial R}{\partial
  u}\,du + R\dphidv\,dv$ with the information from using
  $R\dphidu\,du + \phi\frac{\partial R}{\partial v}\,dv
  $. For our goal it will be better not to symmetrize in
  this manner. Let us recall as a warming-up how one can use
  Riemann's method to find $\phi(t,x)$ for $t>0$ when $\phi$
  and $\dphidt$ are known for $t=0$. Let $P = (t,x)$, $A =
  (0,x-t)$, $B = (0,x+t)$, and $R(Q) = R(P-Q)$. 
\[ \phi(P) - \phi(A) = \int_{A\to P} \dphidv dv = \int_{A\to
  P} R \dphidv dv + \phi \frac{\partial R}{\partial
  u} du = \int_{A\to B} + \int_{B\to P} = \int_{A\to B}\]
Hence:
\[ \phi(P) = \phi(A) + \int_{A\to B} ( R \dphidv
  + \phi \frac{\partial R}{\partial u})  \frac{dx}2 \]
Using $R\dphidu\,du + \phi\frac{\partial R}{\partial v}\,dv$ we get
in the same manner:
\begin{equation}\label{eq:2} \phi(P) = \phi(B) - \int_{A\to B}  ( R \dphidu
  + \phi \frac{\partial R}{\partial v})  \frac{dx}2 \end{equation}
After averaging:
\[ \phi(P) = \frac{\phi(A) + \phi(B)}2 + \frac12 \int_{A\to
  B} ( R \dphidt - \phi \frac{\partial R}{\partial t}) dx\]
This gives the classical formula ($t>0$):
\begin{equation}\label{eq:classic}
\begin{split} \phi(t,x) = \frac{\phi(0,x-t) +
  \phi(0,x+t)}2 
- \frac12 \int_{x-t}^{x+t} t\; \frac{J_1(\sqrt{t^2 -
  (x-x')^2})}{\sqrt{t^2 - (x-x')^2}}\phi(0,x')dx'\\
+ \frac12 \int_{x-t}^{x+t} J_0(\sqrt{t^2 -
  (x-x')^2})\dphidt(0,x')dx'
\end{split}\end{equation}
I have not tried to use it to establish theorem
  \ref{theo:1}.  Anyway, when $\phi$,
  $\dphidx$, $\dphidt$ all belong to $L^2$ at $t=0$, this
  formula shows that $\phi(P)$ is continuous in $P$
  for $t>0$. Replacing $t=0$ with $t=-T$, we find that
  $\phi$ is continuous on spacetime. 

Let us now consider the problem, with the notations of
  Theorem \ref{theo:4}, of determining $k(v) = \psi(v,v)$
  for $v>0$
  when $F(x) = \phi(0,x) = -\dpsidu(0,x)$ and $G(x) =
  \psi(0,x) = - \dphidv(0,x)$ are known for $x>0$. We use $P
  = (v_0,v_0)$, $A = (0,0)$, $B= (0,2v_0)$. We then have:
\[ R(t,x) = J_0(\sqrt{(v_0-t)^2 - (v_0-x)^2}) =
  J_0(2\sqrt{u(v_0-v)})\qquad R(0,x) = J_0(\sqrt{x(2v_0 - x)})\]
\[ \frac{\partial R}{\partial v} =
  \frac{J_1(2\sqrt{u(v_0-v)})}{2\sqrt{u(v_0-v)}} 2u \qquad
  \frac{\partial R}{\partial v}(0,x) =
  \frac{J_1(\sqrt{x(2v_0 - x)})}{\sqrt{x(2v_0 - x)}} x\]

Hence, using \eqref{eq:2} (for $\psi$):
\begin{equation} \psi(v,v) = G(2v) + \frac12\int_0^{2v} ( 
  J_0(\sqrt{x(2v_0 - x)}) F(x) - x\, \frac{J_1(\sqrt{x(2v_0
  - x)})}{\sqrt{x(2v_0 - x)}} G(x)) \,dx\end{equation}

We then consider the converse problem of expressing $G(x) =
  \psi(0,x)$ in terms of $k(v) = \psi(v,v)$. We choose
  $x_0>0$, and consider the rectangle with vertices $P =
  (\frac12 x_0, \frac12 x_0)$, $Q = (0, x_0)$, $Q' = (X,
  X+x_0)$, $P' = (X+ \frac12x_0, X+\frac12 x_0)$ for
  $X\gg0$. We take Riemann's function $S$ to be $1$ on the
  edges $P\to Q$ and $Q\to Q'$. We then write:
\[ \psi(Q) - \psi(P) = \int_{P\to Q}\dpsidu \,du =
  \int_{P\to Q}(S \dpsidu \,du + \psi \frac{\partial
  S}{\partial v} dv)\]
\[ = \int_{P\to P'} + \int_{P'\to Q'} + \int_{Q'\to Q} =
  \int_{P\to P'} \psi  \frac{\partial
  S}{\partial v}\, dv  + \int_{P'\to Q'} S\dpsidu\,du\]
\begin{equation} G(x_0) = \psi(\frac{x_0}2,\frac{x_0}2) + \int_{P\to P'}
  \psi   \frac{\partial 
  S}{\partial v}\, dv - \int_{P'\to Q'} S\phi \,du\end{equation}
Now, $|S|\leq 1$ on the segment leading from $P'$ to $Q'$,
  so we can bound the last integral, using Cauchy-Schwarz,
  then the energy integral, and finally the theorem
  \ref{theo:1}. So this term goes to $0$. On the light cone
  half line from $P$ to $\infty$
  we have:
\[ S(v,v) = J_0(\sqrt{x_0(2v-x_0)})\qquad  \frac{\partial
  S}{\partial v} = -
  \frac{J_1(\sqrt{x_0(2v-x_0)})}{\sqrt{x_0(2v-x_0)}}\,x_0\]
\begin{equation} G(x_0) = \psi(\frac{x_0}2,\frac{x_0}2) -
  \int_{x_0/2}^\infty
  \frac{J_1(\sqrt{x_0(2v-x_0)})}{\sqrt{x_0(2v-x_0)}}\,x_0
  \psi(v,v)\,dv
\end{equation}

Our last task is to obtain the formula for $F(x_0)$. We use
  the same rectangle and same function $S$.
\[ \phi(Q) - \phi(Q') = \int_{Q'\to Q} \dphidv dv =
  \int_{Q'\to Q}  S\dphidv dv +  \phi \frac{\partial
  S}{\partial u} du = \int_{Q'\to P'} \phi  \frac{\partial
  S}{\partial u} du + \int_{P'\to P} S\dphidv dv + 0 \]
On the segment
  $Q'\to P'$ we integrate by parts to get:
\[ \int_{Q'\to P'} \phi  \frac{\partial
  S}{\partial u} du = \phi(P')S(P') - \phi(Q') - \int_{Q'\to
  P'} \dphidu S\,du\]
Again we can bound $S$ by $1$ and apply Cauchy-Schwarz to
$\int_{Q'\to
  P'} \dphidu S\,du$. Then we observe that $\int_{Q'\to P'}
  |\dphidu|^2 |du|$ is bounded above by the energy integral,
  which itself is bounded above by the energy integral on
  the horizontal line having $P'$ as its left end. By
  Theorem \ref{theo:1} this goes to $0$.  And regarding
  $\phi(P')$ one has $\lim_{v\to+\infty} \phi(v,v)=0$ as
  $\phi(v,v)$ and its derivative belong to
  $L^2(0,+\infty;\,dv)$. We cancel the $\phi(Q')$'s on both
  sides of our equations and obtain:
\[ \phi(Q) = - \int_{P\to(\infty,\infty)}  S \dphidv dv =
  +\int_{P\to(\infty,\infty)} S \psi dv \]
Hence
\begin{equation} F(x_0) = 
  \int_{x_0/2}^\infty J_0(\sqrt{x_0(2v-x_0)})\psi(v,v)\,dv\end{equation}

In conclusion: the functions $F(x) = \phi(0,x)$, $G(x) =
  \psi(0,x)$, and $k(v) = \psi(v,v)$ of Theorem
  \ref{theo:4} are related by the following formulas:
  \begin{subequations}
\begin{align}
    F(x) &= 
  \int_{x/2}^\infty J_0(\sqrt{x(2v-x)})k(v)\,dv\\
G(x) &= k(\frac{x}2) -
  \int_{x/2}^\infty x\,
  \frac{J_1(\sqrt{x(2v-x)})}{\sqrt{x(2v-x)}}\,
  k(v)\,dv\\
k(v)  &= G(2v) + \frac12\int_0^{2v} 
  J_0(\sqrt{x(2v - x)}) F(x) \,dx  - \frac12\int_0^{2v}
  x\, \frac{J_1(\sqrt{x(2v 
  - x)})}{\sqrt{x(2v - x)}} G(x) \,dx
\end{align}
  \end{subequations}
Exchanging $F$ and $G$ is like applying a time reversal so
  it corresponds exactly to exchanging $k(v) = \psi(v,v)$
  with $g(u) = \phi(-u,u)$. So the proof of Theorem
  \ref{theo:4} is complete.

\section{Conformal coordinates and concluding remarks}

The Rindler coordinates $(\xi,\eta)$ in the right wedge are defined by
  the equations $x =
  \eta \cosh\xi$, $t = \eta \sinh\xi$. Let us use the
  conformal coordinate system:
  \[
    \xi = \frac12 \log\frac{x+t}{x-t}\qquad\qquad \zeta =
  \frac12 \log(x^2 - t^2) - \log2  = \log\frac\eta2
  \]
where $-\infty<\xi<+\infty$, $-\infty<\zeta<+\infty$. The
variable $\xi$
  plays the r\^ole of time for our scattering. The reason
  for $-\log2$ in $\zeta$ is the following: at $t=0$ this
  gives $e^\zeta = \frac12 x = u = v$. The differential
  equations we shall write are related to the understanding
  of the vanishing condition for an $\cH$ pair on an
  interval $(0,a)$. And $a = \frac12 (2a)$ hence the
  $-\log2$ (to have equations identical with those in
  \cite{hankel}.) The Klein-Gordon equation becomes:
  \begin{equation}
    \label{eq:6}
    \frac{\partial^2\phi}{\partial\xi^2} -
  \frac{\partial^2\phi}{\partial\zeta^2} + 4e^{2\zeta} \phi
  = 0
  \end{equation}
If we now look for ``eigenfunctions'', oscillating
  harmonically in time, $\phi = e^{-i\gamma\xi}\, \Phi(\zeta)$,
  $\gamma\in\RR$, we obtain a Schr\"odinger eigenvalue
  equation:
  \begin{equation}
    \label{eq:4}
     - \Phi''(\zeta) +   4e^{2\zeta} \Phi(\zeta) = \gamma^2 \Phi(\zeta)
  \end{equation}
This Schrödinger operator has a potential function which can
  be conceived of as acting as a repulsive exponential
  barrier for the de~Broglie wave function of a quantum
  mechanical particle coming from $-\infty$ and being
  ultimately bounced back to $-\infty$. The solutions of
  \eqref{eq:4} are the modified Bessel functions
  (\cite{watson}) of imaginary argument $i\gamma$ in the
  variable $2e^{\zeta}$. For each $\gamma\in\CC$ the unique
  (up to a constant factor) solution
 of \eqref{eq:4} which is square integrable at $+\infty$ is
  $K_{i\gamma}(2e^{\zeta})$. 

\From\ Theorem
  \ref{theo:4} it is more convenient to express the $\cH$
  transform as a scattering for the two-component,
  ``Dirac'', differential system. The spinorial nature of
  $\left[\begin{smallmatrix} \psi\\
  \phi\end{smallmatrix}\right]$ leads under the change of
  coordinates $(t,x)\mapsto (\xi,\zeta)$ to
  $e^{\frac\zeta2}e^{-\frac\xi2}\phi$ rather than $\phi$,
  and to $e^{\frac\zeta2}e^{+\frac\xi2}\psi$ rather than
  $\psi$. In order to get quantities which, in the past at
  $\xi\to-\infty$, look like $\phi$ and, in the future at
  $\xi\to+\infty$, look like $\psi$ we consider the linear
  combinations:
  \begin{subequations}
    \begin{align}\label{eq:ABdefa}
      \cA &= \frac12 \,e^{\frac\zeta2}\,(
      +e^{-\frac\xi2}\phi + e^{\frac\xi2} \psi)\\
\label{eq:ABdefb}
      \cB &= \frac i2 \,e^{\frac\zeta2}\,(
      -e^{-\frac\xi2}\phi + e^{\frac\xi2} \psi)
    \end{align}
  \end{subequations}
Their differential system is:
\begin{subequations}
  \begin{align}
    +i\frac{\partial\cA}{\partial\xi} &=
  +\left(\frac\partial{\partial\zeta} - 2e^\zeta\right)\cB
 \\ +i\frac{\partial\cB}{\partial\xi} &= 
  -\left(\frac\partial{\partial\zeta} + 2e^\zeta\right)\cA
  \end{align}
\end{subequations}
Or, if we look for solutions oscillating in time as $e^{-i\gamma\xi}$:
\begin{subequations}\label{eq:ABsys}
  \begin{align}\label{eq:ABsysa}
  \left(\frac\partial{\partial\zeta} - 2e^\zeta\right)\cB &= \gamma\cA
\\ \label{eq:ABsysb}
  \left(-\frac\partial{\partial\zeta} - 2e^\zeta\right)\cA &= \gamma\cB
  \end{align}
\end{subequations}
and this gives Schrödinger equations:
\begin{subequations}
  \begin{align}
\label{eq:ABschroda}
    -  \frac{\partial^2\cA}{\partial\zeta^2} + (4e^{2\zeta} -
  2e^\zeta) \cA  &= \gamma^2 \cA \\ 
\label{eq:ABschrodb}
    -  \frac{\partial^2\cB}{\partial\zeta^2} + (4e^{2\zeta} +
  2e^\zeta) \cB  &= \gamma^2 \cB
   \end{align}
\end{subequations}
So we have two exponential barriers, and two associated
``scattering functions''
  giving the induced phase shifts.  From our previous
  discussion of the scattering in the Lax-Phillips formalism
  we can expect from equation \eqref{eq:chi} that a
  formalism of Jost functions will confirm these functions to be
  \begin{equation}
    \label{eq:S}
    \cal S(\gamma) = \frac{\Gamma(\frac12 -
  i\gamma)}{\Gamma(\frac12 + i \gamma)}\qquad
  (\gamma\in\RR)\;,
  \end{equation}
for the equation associated with $\cA$ and $-\cal
  S(\gamma)$ for the equation associated with $\cB$. And
  indeed the solution $\left[\begin{smallmatrix} \cA_\gamma\\
  \cB_\gamma\end{smallmatrix}\right]$ of the system
  \eqref{eq:ABsys} which is square-integrable at $+\infty$
  is given by the formula
  \begin{equation}
    \label{eq:10}
    \begin{bmatrix} \cA_\gamma(\zeta)\\
  \cB_\gamma(\zeta)\end{bmatrix} = \begin{bmatrix}
  e^{\frac\zeta2}\left(K_s(2e^{\zeta}) +
  K_{1-s}(2e^\zeta)\right)\\ i\,
  e^{\frac\zeta2}\left(K_s(2e^{\zeta}) -
  K_{1-s}(2e^\zeta)\right)\end{bmatrix}\qquad(s=\frac12+i\gamma)
  \end{equation}
Let $j_\gamma(\zeta)$ be the
  solution of \eqref{eq:ABschroda} which satisfies the Jost
  condition  $j_\gamma(\zeta)\sim
  e^{-i\gamma\zeta}$ as $\zeta\to-\infty$. Then the exact
  relation holds (a detailed
  treatment is given in \cite{hankel}):
  \begin{equation}
    \label{eq:8}
    \cA_\gamma(\zeta) = \frac12\left(\Gamma(s) j_\gamma(\zeta)
  + \Gamma(1-s) j_{-\gamma}(\zeta)\right)\qquad(s=\frac12+i\gamma)
  \end{equation}
We interpret this as saying that the $\cA$-wave comes
  from $-\infty$ and is bounced back with a phase-shift
  which at frequency $\gamma$ equals
  $\arg\frac{\Gamma(\frac12
  -i\gamma)}{\Gamma(\frac12+i\gamma)} = \arg\cS(\gamma)
  $. For the $\cB$ equation one obtains $-\cS(\gamma)$ as
  the phase shift function.

We have
  associated in \cite{cras4} Schr\"odinger equations to the
  cosine and sine kernels whose potential functions also
  have exponential vanishing at $-\infty$ and exponential
  increase at $+\infty$, and whose associated scattering
  functions are the functions
  arising in the functional equations of the Riemann and
  Dirichlet L-functions. The equations (13a), (13b) of
  \cite{cras4} are analogous to \eqref{eq:ABdefa},
  \eqref{eq:ABdefb} above, and (14a), (14b) of \cite{cras4}
  are analogous to \eqref{eq:ABsysa} and \eqref{eq:ABsysb}
  above. The analogy is no accident. The reasoning of
  \cite{cras4} leading to the consideration of Fredholm
  determinants when trying to understand self- and
  skew-reciprocal functions under a scale reversing operator
  on $L^2(0,+\infty;\,dx)$ is quite general. The (very
  simple) potential functions in the equations
  \eqref{eq:ABschroda} and \eqref{eq:ABschrodb} can be
  written  in terms of Fredholm determinants associated
  with the $\cH$ transform. The detailed treatment is given in
  \cite{hankel}.

The function $\cal S(\gamma)$ arises in
number theoretical
  functional equations (for the Dedekind zeta functions of
  imaginary quadratic fields). We don't know if its
  interpretation obtained here in terms of the Klein-Gordon
  equation may lead us to legitimately hope for number
  theoretical applications. An interesting physical context
  where $S(\gamma)$ has appeared is the method of angular
  quantization in integrable quantum field theory
  \cite[App. B]{lukzak}. And, of course the group
  of Lorentz boosts and the Rindler wedge are
  connected by the Bisognano-Wichman theorem
  \cite{biswich1,biswich2,haag}.

The potentials associated in \cite{cras4} to the cosine and
  sine kernels are, contrarily to the simple-minded
  potentials obtained here, mainly known through their
  expressions as Fredholm determinants, and these are
  intimately related to the Fredholm determinant of the
  Dirichlet kernel, which has been found to be so important
  in random matrix theory.
It is thus legitimately considered an important problem to
try to acquire for
  the cosine and sine kernels the kind of understanding
  which has been achieved here for the $\cH$ transform. Will
  it prove possible to achieve this on (a subset, with
  suitable conformal coordinates) of (possibly higher
  dimensional) Minkowski space?

We feel that some kind of
  non-linearity should be at work. 
A tantalizing thought presents itself: perhaps the
  kind of understanding of the Fourier transform which is
  hoped for will arise from the study of the causal
  propagation and scattering of (quantum mechanical?) waves
  on a certain curved Einsteinian spacetime.

\setlength{\parskip}{\the\baselineskip}
\addtolength{\baselineskip}{-6pt}

\end{document}